\documentclass{amsart}
\usepackage{amsmath}
\usepackage{amsfonts}
\usepackage{amssymb}

\newtheorem{theorem}{Theorem}
\newtheorem*{lemma}{Lemma}
\newtheorem*{proposition}{Proposition}

\begin{document}

\title{The Lie algebra of cyclic coinvariants of a symplectic space}

\author[E.~Kushnirsky]{Eugene Kushnirsky}
\address{Department of Mathematics,
Northwestern University,
2033 Sheridan Road,
Evanston, IL 60208}
\email{ekushnir@math.northwestern.edu}

\author[M.~Larsen]{Michael Larsen}
\address{Department of Mathematics, Indiana University,
Bloomington, IN 47405}
\email{larsen@math.indiana.edu}

\begin{abstract} We exhibit a natural Lie algebra structure on the graded space
  of cyclic coinvariants of a symplectic vector space.

\end{abstract}


\maketitle

\section{Introduction}
Let $V$ be a vector space endowed with a symplectic form.  It is
well known that the adjoint representation of $\textrm{Sp}(V)$ is
equivalent 
to the symmetric square of $V$, i.e., to the coinvariants of 
$\mathbb Z/2\mathbb Z$ acting by cyclic permutation on $V\otimes V$.
This allows us 
to endow $(V^{\otimes 2})_{\mathbb Z/2\mathbb Z}$ with the structure
of a (symplectic) Lie 
algebra.  In this note, we show that this is just the degree $0$ part of
a larger object: the graded vector space
$$\bigoplus_{n=0}^\infty(V^{\otimes n+2})_{\mathbb Z/(n+2)\mathbb Z}$$
has a natural structure
of Lie algebra  
determined by the symplectic form on $V$.   In the special case that
$\dim V = 2$, for every associative algebra $A$, this algebra acts
infinitesimally on fibers of the commutator map $[\;,\;]:\;A^2\to A$.

This paper developed from our interest in the fibers of the commutator
map on the Lie algebra $\mathfrak{sl}_n$.  After we wrote it,
B.~Tsygan called our attention to a paper of M.~Kontsevich \cite{K}
which arrived at similar conclusions, coming from a rather different
(and more sophisticated) point of view.  Kontsevich considered the
tensor algebra of a symplectic space $V$ over $\mathbb Q$ and looked
at the Lie algebra of derivations preserving its symplectic form.  He
asserted that this algebra is naturally isomorphic to
$$\bigoplus_{n=0}^\infty(V^{\otimes n+2})^{{\mathbb Z}/(n+2){\mathbb
    Z}}.$$ 
Of course in characteristic zero, there is a natural isomorphism
  between invariants and coinvariants for any finite group, which is
    not the case in characteristic $p$.  By contrast, we give a
    characteristic free formula for the Lie bracket which is naturally  defined on the \emph{coinvariant} space.  
In addition, we indicate a \emph{dihedral} theory which may be useful
    in analyzing the commutator fibers for $\mathfrak{so}_n$ and $\mathfrak{sp}_{2n}$.

\section{The Lie Algebra}
Let $K$ be a field and $V$ a vector space over $K$ equipped with an
alternating bilinear form $\langle\thinspace,\thinspace\rangle$. For any $l$,
$\mathbb{Z}/l\mathbb{Z}$ acts on $V^{\otimes l}$ in the obvious way and we
form the graded $K$-vector space

\[
L(V)=  \bigoplus_{l=2}^{\infty}(V^{\otimes l})_{\mathbb{Z}/l\mathbb{Z}%
}  \;.
\]

\noindent where $(V^{\otimes l})_{\mathbb{Z}/l\mathbb{Z}}$ has degree $l-2$.
We now define a Lie bracket on $L(V)$. 
To do this, let us first
introduce some additional notation.  Fix a basis $\mathcal B$ of $V$.
Given integers $i$ and $j$ such that
$1\leq i,j\leq l$, we set 
$$D_{i,j}(\alpha_{1}\otimes\cdots\otimes\alpha_{l})
=\begin{cases}
\alpha_{i+1}\otimes\cdots\otimes\alpha_{j-1}&\text{if $i+1\le j-1$,}\\
0&\text{if $i = j-1$},\\
\alpha_{i+1}\otimes\cdots\otimes\alpha_{l}\otimes\alpha_{1}\otimes
\cdots\otimes\alpha_{j-1}&\text{if $i\ge j$},\\
\end{cases}
$$
where $\alpha_i\in\mathcal B$.  Now 
if $\alpha=\alpha_{1}\otimes\cdots
\otimes\alpha_{l}$ and $\beta=\beta_{1}\otimes\cdots\otimes\beta_{m}$
are tensor products of basis elements, we define
\[
{\lbrack\alpha,\beta]=
\sum_{i=1}^{l}\sum_{j=1}^{m}\langle\alpha_{i},\beta_{j}\rangle\overline{D_{i,i}
(\alpha)\otimes D_{j,j}(\beta)}\;.}
\]
and extend by linearity.  (Here $\bar\gamma$ denotes the class of
$\gamma\in V^{\otimes n}$ in $(V^{\otimes n})_{\mathbb{Z}/n\mathbb{Z}}$.)
Note that if $\sigma$ denotes the canonical generator of
$\mathbb{Z}/l\mathbb{Z}$ then 
\[
D_{i,i}(\sigma^r(\alpha))=
\begin{cases}
D_{r+i,r+i}(\alpha) & \text{if $i\leq l-r$}\\
D_{r+i-l,r+i-l}(\alpha) & \text{if $i>l-r$.}
\end{cases}
\]
Thus
$[\;,\;]:(V^{\otimes l})_{\mathbb{Z}/l\mathbb{Z}}
\times(V^{\otimes m})_{\mathbb{Z}/m\mathbb{Z}}
\longrightarrow(V^{\otimes(l+m-2)})_{\mathbb{Z}/(l+m-2)\mathbb{Z}}$ is well-defined.
Moreover, it does not depend on the choice of basis $\mathcal B$.

For each $x\in\mathcal B$, we also define $D_x:L(V)\to L(V)$ by
\[
D_x(\alpha)=\sum_{i=1}^l\delta_{\alpha_i,x}\overline{D_{i,i}(\alpha)},
\]
where $\delta$ is the Kronecker delta and $\alpha$ is as before. We
extend $D_x$ to a (well-defined) endomorphism of $L(V)$.

\begin{theorem}
The above bracket makes $L(V)$ into a Lie algebra. Moreover each $D_x$
is a derivation.
\end{theorem}

\begin{proof}
It is easily verified that the bracket is bilinear and antisymmetric. We
now prove that the Jacobi identity holds. 
Suppose 
$\alpha=\alpha_{1}\otimes\cdots
\otimes\alpha_{l}$, $\beta=\beta_{1}\otimes\cdots\otimes\beta_{m}$ and 
$\gamma=\gamma_{1}\otimes\cdots\otimes\gamma_{n}$; then we
have
\begin{gather*}
\lbrack\alpha,[\beta,\gamma]]+[\beta,[\gamma,\alpha]]+[\gamma,[\alpha
,\beta]]=\\
\sum_{i=1}^{l}\sum_{j=1}^{m}\sum_{k=1}^{n}\sum_{\substack{r=1\\r\neq j}}^{m}
\langle\beta_{j},\gamma_{k}\rangle\langle\alpha_{i},\beta_{r}\rangle\overline{D_{i,i}
(\alpha)\otimes D_{r,j}(\beta)\otimes D_{k,k}(\gamma)\otimes D_{j,r}(\beta)}\\
+\sum_{i=1}^{l}\sum_{j=1}^{m}\sum_{k=1}^{n}\sum_{\substack{r=1\\r\neq k}}^{n}
\langle\beta_j,\gamma_k\rangle\langle\alpha_{i},\gamma_{r}\rangle\overline{D_{i,i}%
(\alpha)\otimes D_{r,k}(\gamma)\otimes D_{j,j}(\beta)\otimes D_{k,r}(\gamma
)}\\
+\sum_{i=1}^{l}\sum_{j=1}^{m}\sum_{k=1}^{n}\sum_{\substack{r=1\\r\neq k}%
}^{n}\langle\gamma_{k},\alpha_{i}\rangle\langle\beta_{j},\gamma_{r}\rangle\overline{D_{j,j}%
(\beta)\otimes D_{r,k}(\gamma)\otimes D_{i,i}(\alpha)\otimes D_{k,r}(\gamma
)}\\
+\sum_{i=1}^{l}\sum_{j=1}^{m}\sum_{k=1}^{n}\sum_{\substack{r=1\\r\neq i}%
}^{l}\langle\gamma_{k},\alpha_{i}\rangle\langle\beta_{j},\alpha_{r}\rangle\overline{D_{j,j}%
(\beta)\otimes D_{r,i}(\alpha)\otimes D_{k,k}(\gamma)\otimes D_{i,r}(\alpha
)}\\
+\sum_{i=1}^{l}\sum_{j=1}^{m}\sum_{k=1}^{n}\sum_{\substack{r=1\\r\neq i}%
}^{l}\langle\alpha_{i},\beta_{j}\rangle\langle\gamma_{k},\alpha_{r}\rangle\overline{D_{k,k}%
(\gamma)\otimes D_{r,i}(\alpha)\otimes D_{j,j}(\beta)\otimes D_{i,r}(\alpha
)}\\
+\sum_{i=1}^{l}\sum_{j=1}^{m}\sum_{k=1}^{n}\sum_{\substack{r=1\\r\neq j}%
}^{m}\langle\alpha_{i},\beta_{j}\rangle\langle\gamma_{k},\beta_{r}\rangle\overline{D_{k,k}%
(\gamma)\otimes D_{r,j}(\beta)\otimes D_{i,i}(\alpha)\otimes
  D_{j,r}(\beta)}\thinspace .%
\end{gather*}

\noindent Now note that 
$$\overline{D_{i,i}(\alpha)\otimes D_{r,j}(\beta)\otimes
D_{k,k}(\gamma)\otimes D_{j,r}(\beta)}=\overline{D_{k,k}(\gamma)\otimes
D_{j,r}(\beta)\otimes D_{i,i}(\alpha)\otimes D_{r,j}(\beta)}$$
in $(V^{\otimes(l+m+n-4)})_{\mathbb{Z}/(l+m+n-4)\mathbb{Z}}$. Therefore,
interchanging the sums over $j$ and $r$ in the last term, we see that it is
exactly the negative of the first term (since $\langle\alpha_{i},\beta
_{j}\rangle\langle\gamma_{k},\beta_{r}\rangle=
-\langle\beta_{r},\gamma_{k}\rangle\langle\alpha_{i},\beta_{j}\rangle$).   
Similarly, the second term cancels the third and the fourth term
cancels 
the fifth. The Jacobi identity now follows by linearity.

Finally, we have 
\begin{eqnarray*}
D_x(\lbrack\alpha,\beta\rbrack) &=&
\sum_{i=1}^l\sum_{j=1}^m\langle\alpha_{i},\beta_{j}\rangle
D_x(\overline{D_{i,i}(\alpha)\otimes
    D_{j,j}(\beta)})\\
&=& \sum_{i=1}^l\sum_{j=1}^m
\langle\alpha_{i},\beta_{j}\rangle\sum_{\substack{k=1\\k\neq
    i}}^{l}\delta_{\alpha_k,x}\overline{D_{k,i}(\alpha)\otimes
    D_{j,j}(\beta)\otimes
  D_{i,k}(\alpha)}\\
&+& \sum_{i=1}^l\sum_{j=1}^m
\langle\alpha_{i},\beta_{j}\rangle
\sum_{\substack{k=1\\k\neq
    j}}^{m}\delta_{\beta_k,x}\overline{D_{k,j}(\beta)\otimes
    D_{i,i}(\alpha)\otimes 
  D_{j,k}(\beta)}\\
&=& \lbrack D_x(\alpha),\beta\rbrack+ \lbrack\alpha,D_x(\beta)\rbrack 
\end{eqnarray*}
\end{proof}

If $\phi:V\to W$ is a linear transformation such that
$\langle\phi(\alpha),\phi(\beta)\rangle_W =
\langle\alpha,\beta\rangle_V$ for all 
$\alpha,\beta\in V$ then we have an obvious Lie algebra homomorphism
$L(\phi):L(V)\to L(W)$ induced by
$\alpha_{1}\otimes\cdots\otimes\alpha_{l}\mapsto
\phi(\alpha_{1})\otimes\cdots\otimes\phi(\alpha_{l})$. 

\vspace{.25in}
The case $\dim V = 2$ is of particular interest.
Let $\lbrace x,y\rbrace$ be
a basis with $\langle x,y\rangle=1$.  We identify $V^{\otimes l}$ with the 
noncommutative homogeneous polynomials of degree $l$ in $x$ and $y$.

Suppose that $A$ is an associative $K$-algebra 
considered as a Lie algebra via the usual bracket $[X,Y]=XY-YX$.
Let $N_l$ denote the
the trace with
respect to the $\mathbb Z/l\mathbb Z$-action on $V^{\otimes l}$; then
the $N_l$ induce a well-defined map 
$N:L(V)\to\bigoplus_{l\ge 2} (V^{\otimes l})^{\mathbb Z/l\mathbb Z}$. 
Let $K[\varepsilon]$ denote the ring of dual numbers. 
If $\alpha\in V^{\otimes l}$ represents a class in $L(V)$, write
$N(\alpha)=xp_\alpha-yq_\alpha$ and 
consider the vector field 
$$
F_\alpha(X,Y)=(X-\varepsilon q_\alpha(X,Y),Y-\varepsilon p_\alpha(X,Y))
$$
on $A\times A$.  Note that $N(\alpha)$, $p_\alpha$,
$q_\alpha$, 
and $F_\alpha$ depend only on $\bar\alpha$.

\begin{theorem}
The map $\bar\alpha\mapsto F_\alpha$ is a homomorphism of
Lie algebras from $L(V)$ to the algebra of vector fields
on $A\times A$ tangent to the fibers of
the commutator map.
\end{theorem}

\begin{proof}
As $\sigma(N_l(\alpha)) = N_l(\alpha)$, $(xp_\alpha - yq_\alpha) 
= (p_\alpha x-q_\alpha y)$, so
$$[X-\varepsilon q_\alpha(X,Y),Y-\varepsilon p_\alpha(X,Y)]\equiv [X,Y]
\pmod{\varepsilon}.$$ 
Thus, $F_\alpha$ is always tangent to the fibers of the commutator
map.  To see that $\bar\alpha\mapsto F_\alpha$ is a Lie algebra homomorphism,
it suffices to check the case that
$\alpha=\alpha_1\otimes\cdots\otimes\alpha_l$ and  
$\beta=\beta_1\otimes\cdots\otimes\beta_m$ are tensor products of
basis vectors. 
Now,
$$p_\alpha = \sum_{i=1}^l \delta_{\alpha_i,x}D_{i,i}(\alpha),
\ q_\alpha = -\sum_{i=1}^l \delta_{\alpha_i,y}D_{i,i}(\alpha).$$
(Note that $\overline{p_{\alpha}}=D_x(\alpha)$ and
$\overline{q_{\alpha}}=D_y(\alpha)$.) 
Regarding $F_\alpha$ and $F_\beta$ as sections 
$A\times A\to A[[\varepsilon]]\times A[[\varepsilon]]$ of 
the evaluation at zero map,
\begin{align*}
F_\alpha(&F_\beta(X,Y)) 
= F_\alpha\Bigl(X+\varepsilon\sum_{j=1}^m
\delta_{\beta_j,y}D_{j,j}(\beta)(X,Y),  
    Y-\varepsilon\sum_{j=1}^m \delta_{\beta_j,x}D_{j,j}(\beta)(X,Y)\Bigr) = \\
  \Bigl(X&+\varepsilon\sum_{j=1}^m \delta_{\beta_j,y}D_{j,j}(\beta)(X,Y) \\
    +\varepsilon&\sum_{i=1}^l \delta_{\alpha_i,y}D_{i,i}(\alpha)
    \Bigl(X+\varepsilon\sum_{j=1}^m \delta_{\beta_j,y}D_{j,j}(\beta)(X,Y), 
    Y-\varepsilon\sum_{j=1}^m \delta_{\beta_j,x}D_{j,j}(\beta)(X,Y)\Bigr), \\
Y&-\varepsilon\sum_{j=1}^m \delta_{\beta_j,x}D_{j,j}(\beta)(X,Y)   \\  
    -\varepsilon&\sum_{i=1}^l \delta_{\alpha_i,x}D_{i,i}(\alpha)
    \Bigl(X+\varepsilon\sum_{j=1}^m \delta_{\beta_j,y}D_{j,j}(\beta)(X,Y), 
    Y-\varepsilon\sum_{j=1}^m
    \delta_{\beta_j,x}D_{j,j}(\beta)(X,Y)\Bigr)\Bigr). \\ 
\end{align*}

If $\gamma=\gamma_1\otimes\cdots\otimes\gamma_n$ 
is a tensor monomial regarded as a noncommutative
homogeneous polynomial, then
\begin{multline*}
\gamma(X+\varepsilon X_1,Y+\varepsilon Y_1)\equiv\\
\gamma(X,Y) + \varepsilon\sum_{k=1}^n D_{n,k}(\gamma)(X,Y)
( \delta_{\gamma_k,x} X_1 + \delta_{\gamma_k,y}
Y_1)D_{k,1}(\gamma)(X,Y)\pmod{\varepsilon^2}.\end{multline*} 
Therefore,

\begin{align*}
F_\alpha(F_\beta & (X,Y)) -  F_\beta(F_\alpha(X,Y))\\
    \equiv \varepsilon^2 \Bigl(
 &\sum_{i=1}^l \sum_{j=1}^m\sum_{\substack{k=1\\k\neq i}}^l
    \delta_{\alpha_i,y} 
    (\delta_{\alpha _k,x}\delta_{\beta_j,y} -  \delta_{\alpha
    _k,y}\delta_{\beta_j,x}) 
    \bigl(D_{i,k}(\alpha)D_{j,j}(\beta)D_{k,i}(\alpha)\bigr)(X,Y),\\
& \sum_{i=1}^l \sum_{j=1}^m\sum_{\substack{k=1\\k\neq i}}^l\delta_{\alpha_i,x}
    (-\delta_{\alpha _k,x}\delta_{\beta_j,y} + \delta_{\alpha
    _k,y}\delta_{\beta_j,x}) 
    \bigl(D_{i,k}(\alpha)D_{j,j}(\beta)D_{k,i}(\alpha)\bigr)(X,Y)\Bigr)\\
   - \varepsilon^2\Bigl(
& \sum_{j=1}^m \sum_{i=1}^l \sum_{\substack{k=1\\k\neq j}}^m
    \delta_{\beta_j,y} 
    (\delta_{\beta _k,x}\delta_{\alpha_i,y} - \delta_{\beta
    _k,y}\delta_{\alpha_i,x}) 
    \bigl(D_{j,k}(\beta)D_{i,i}(\alpha)D_{k,j}(\beta)\bigr)(X,Y),\\
& \sum_{j=1}^m \sum_{i=1}^l \sum_{\substack{k=1\\k\neq j}}^m
    \delta_{\beta_j,x} 
    (-\delta_{\beta _k,x}\delta_{\alpha_i,y} + \delta_{\beta
    _k,y}\delta_{\alpha_i,x}) 
    \bigl(D_{j,k}(\beta)D_{i,i}(\alpha)D_{k,j}(\beta)\bigr)(X,Y)\Bigr)
\end{align*}

\begin{align*}
\equiv  \varepsilon^2\Bigl(
    &\sum_{i=1}^l \sum_{j=1}^m\sum_{\substack{k=1\\k\neq i}}^l 
     \delta_{\alpha_i,y}
    \langle\alpha _k,\beta_j\rangle
    \bigl( D_{i,k}(\alpha)D_{j,j}(\beta)D_{k,i}(\alpha)\bigr)(X,Y)\\ 
-&    \sum_{j=1}^m \sum_{i=1}^l \sum_{\substack{k=1\\k\neq j}}^m
     \delta_{\beta_j,y}
    \langle\beta_k,\alpha_i\rangle
    \bigl( D_{j,k}(\beta)D_{i,i}(\alpha)D_{k,j}(\beta)\bigr)(X,Y),\\ 
&\sum_{i=1}^l \sum_{j=1}^m\sum_{\substack{k=1\\k\neq i}}^l 
     \delta_{\alpha_i,x}
    \langle\beta_j,\alpha _k\rangle
    \bigl( D_{i,k}(\alpha)D_{j,j}(\beta)D_{k,i}(\alpha)\bigr)(X,Y)\\ 
- &   \sum_{j=1}^m \sum_{i=1}^l \sum_{\substack{k=1\\k\neq j}}^m
     \delta_{\beta_j,x}
    \langle\alpha_i,\beta_k\rangle
    \bigl( D_{j,k}(\beta)D_{i,i}(\alpha)D_{k,j}(\beta)\bigr)(X,Y)\Bigr)\\
 \equiv \varepsilon^2 & (-q_{_{[\bar\alpha,\bar\beta]}}
    (X,Y),-p_{_{[\bar\alpha,\bar\beta]}}(X,Y))\qquad\pmod{\varepsilon^3}.\\
\end{align*}
\end{proof}
\vspace{.25in}

Each $V^{\otimes l}$ has a natural involution $\iota$ induced by
$$\iota(\alpha_{1}\otimes\cdots \otimes\alpha_{l})=
(-1)^{l}\alpha_{l}\otimes\cdots 
\otimes\alpha_{1}.$$
Extending $\iota$ to all of $L(V)$, we get a Lie
algebra involution; to see this, observe that 

\begin{align*}
\iota(\overline{D_{i,i}
(\alpha)\otimes D_{j,j}(\beta)}) & = (-1)^{l+m}\overline{D_{m-j+1,m-j+1}
(\iota(\beta))\otimes D_{l-i+1,l-i+1}(\iota(\alpha))}\\[.1in]
& = (-1)^{l+m}\overline{D_{l-i+1,l-i+1}
(\iota(\alpha))\otimes D_{m-j+1,m-j+1}(\iota(\beta))}\\
\end{align*}
and since $(\iota(\alpha))_{l-i+1}=\alpha_i$,
$(\iota(\beta))_{m-j+1}=\beta_j$ and $(-1)^{l+m-2}=(-1)^{l+m}$, it
follows that
$[\iota(\alpha),\iota(\beta)]=\iota([\alpha,\beta])$. We will also
need to consider the unsigned involution $I$ induced by
$I(\alpha_{1}\otimes\cdots \otimes\alpha_{l}) =
\alpha_{l}\otimes\cdots
\otimes\alpha_{1}$.

Denote by $P_+(V)$ and $P_-(V)$, respectively, the $+1$ and $-1$ eigenspaces
for $\iota$. Then of course $P_+(V)$ is a Lie subalgebra of $L(V)$ and
$[P_-(V),P_-(V)]\subset P_+(V)$. 

Similarly, we let $P_{l,+}$ and $P_{l,-}$ denote the $+1$ and $-1$
eigenspaces for $\iota$ in $V^{\otimes l}$. To simplify the notation,
we will assume that $V$ is finite-dimensional, so that ${\mathcal B}=
\lbrace x_1,\dots, x_r\rbrace$.

\begin{lemma}
Suppose $\alpha\in (V^{\otimes l})^{\mathbb Z/l\mathbb Z}$ and write
$\alpha=\sum x_i\otimes p_i$, where the $p_i\in V^{\otimes l-1}$. Then
$\alpha\in P_{l,+}$ (resp., $P_{l,-}$) if and only if $p_1,\dots,
p_r\in P_{l-1,-}$ (resp., $P_{l-1,+})$.
\end{lemma}

\begin{proof}
This is clear, as $\iota(\alpha)=-\sum\iota(p_i)\otimes x_i$ and
$\sigma(\alpha)= \sum p_i\otimes x_i$.
\end{proof}

\vspace{.25in}
As an application of these definitions, we turn to the case of
 a Lie algebra
$\mathfrak g$ defined over $K$ and identify $V^{\otimes l}$ with
noncommutative polynomials of degree $l$ in $x_1,\dots,x_r$. 
Our first
observation is that if $\mathfrak g$ is either $\mathfrak{so}_n$ or
$\mathfrak{sp}_n$ and $p=p(x_1,\dots,x_r)$ is in 
$\bigoplus_lP_{l,-}$ then for all $X_1,\dots,
X_r\in\mathfrak g$, $p(X_1,\dots,X_r)\in\mathfrak g$. 
Indeed write $p=\sum p_i$, where each $p_i$ is a
monomial. Then if $\mathfrak g=\mathfrak{so}_n$,

\begin{align*}
p(X_1,\dots,X_r)^t & =\sum p_i(X_1,\dots,X_r)^t = \sum
I(p_i)(X_1^t,\dots,X_r^t)\\
& = \sum (-1)^{\deg p_i}I(p_i)(X_1,\dots,X_r) =
\iota(p)(X_1,\dots,X_r)\\
& = -p(X_1,\dots,X_r).\\
\end{align*}
If $\mathfrak g=\mathfrak{sp}_n$ let $J$ denote the matrix of the
nondegenerate alternating form which defines $\mathfrak g$. Then

\begin{align*}
Jp(X_1,\dots,X_r) & =\sum Jp_i(X_1,\dots,X_r) 
 = \sum (-1)^{\deg p_i}p_i(X_1^t,\dots,X_r^t)J\\
& =(\iota(p)(X_1,\dots,X_r))^tJ
 = -p(X_1,\dots,X_r)^tJ.\\
\end{align*}

With $\dim V = 2$ and
$\mathfrak g$ of the above type (in particular, $\mathfrak
g\subset\mathfrak{gl}_n$) we consider the above map $\bar\alpha\mapsto
F_{\alpha}$ from $L(V)$ to the algebra of vector fields on
$\mathfrak{gl}_n\times\mathfrak{gl}_n$.

\begin{proposition}
The image of $P_+(V)$ under this map consists of vector fields on 
$\mathfrak g\times\mathfrak g$ tangent to the fibers of the
commutator map $\mathfrak g\times\mathfrak g\to\mathfrak g$.
\end{proposition}

\begin{proof}
If $\alpha\in V^{\otimes l}$
then $\iota(\sigma(\alpha))=\sigma^{-1}(\iota(\alpha))$, which implies
that $N(P_+(V))\subset\bigoplus_lP_{l,+}$. Moreover if $\alpha\in P_+(V)$
and
$N(\alpha)=xp_{\alpha}-yq_{\alpha}$ then
$p_{\alpha},q_{\alpha}\in\bigoplus_lP_{l,-}$ by the 
Lemma. Now we are done, as
$p_{\alpha}(X,Y),q_{\alpha}(X,Y)\in\mathfrak g$ whenever
$X,Y\in\mathfrak g$ by the preceding discussion.
\end{proof}

\end{document}